\documentclass[10pt,twoside]{article}
\usepackage{amssymb}
\usepackage{setspace}
\usepackage{latexsym}
\usepackage{amsthm}
\usepackage{amsmath}
\usepackage{colortbl}
\usepackage{latexsym}
\newtheorem{theorem}{Theorem}[section]
\newtheorem{definition}{Definition}[section]
\newtheorem{lemma}{Lemma}[section]
\newtheorem{proposition}{Proposition}[section]
\newtheorem{remark}{Remark}[section]

\newtheorem{corollary}{Corollary}[section]
\def\proof{\mbox {\it Proof.~}}
\makeatletter\def\theequation{\arabic{section}.\arabic{equation}}\makeatother
\newcommand{\be} {\begin{equation}}
\newcommand{\ee} {\end{equation}}
\newcommand{\bea} {\begin{eqnarray}}
\newcommand{\eea} {\end{eqnarray}}
\newcommand{\Bea} {\begin{eqnarray*}}
\newcommand{\Eea} {\end{eqnarray*}}

\newcommand{\e} {\epsilon}

\newcommand{\al} {\alpha}

\newcommand{\de} {\delta}
\newcommand{\ga} {\gamma}

\newcommand{\Om} {\Omega}

\newcommand{\De} {\Delta}
\newcommand{\la} {\lambda}

\newcommand{\La} {\Lambda}

\newcommand{\noi} {\noindent}

\newcommand{\va} {\varphi}
\newcommand{\rar} {\rightarrow}

\newcommand{\HH}{H^1_0(\Om)}
\newcommand{\II}{I_\la}

\newcommand{\RR}{{\mathbb R}}

\makeatletter\def\theequation{\arabic{section}.\arabic{equation}}\makeatother
\begin{document}
\title{
{\bf\Large  Elliptic Problems in $\mathbb{R}^N$
 with  Critical and Singular Discontinuous Nonlinearities}}
\author{{\bf\large R. Dhanya}\footnote{The author acknowledges the support of National Board of Higher Mathematics,
DAE , Govt. of India for providing financial support under the grant no. R(IA)-NBHM-PDF(DR-MA)/2013-4126.}\\
{\it\small Department of Mathematics, IISc, Bangalore 560012}\\
{\it\small e-mail: dhanya.tr@gmail.com}\\
\vspace{1mm}
{\bf\large S. Prashanth} and
{\bf\large Sweta Tiwari}\vspace{1mm}\\
{\it\small TIFR-CAM, Post Bag No.6503, Sharada
Nagar,}\\ {\it\small Chikkabommasandra, Bangalore 560065.}\\
{\it\small e-mail:pras@math.tifrbng.res.in, sweta@math.tifrbng.res.in}\\
\vspace{1mm}
{\bf\large K. Sreenadh}\\
{\it\small Department of Mathematics, Indian Institute of Technology Delhi}\\ {\it \small Hauz Khas, New Delhi 110016.}\\
{\it\small e-mail: sreenadh@gmail.com}}
\date{}
\maketitle
\begin{center}
{\bf\small Abstract}
\vspace{3mm}
\hspace{.05in}\parbox{4.5in}
{{\small \noindent Let $\Omega$ be a bounded domain in $\mathbb R^{N}$, $N\geq3$ with smooth boundary,
 $a>0, \lambda>0$ and $0<\delta<3$ be real numbers. Define $2^*:=\displaystyle\frac{2N}{N-2}$ and the characteristic function of a set $A$ by $\chi_A$.
 We consider the following critical
problem with singular and discontinuous nonlinearity:
\begin{eqnarray*}
 (P_\la^a)~~~~ \qquad \Biggl\{\begin{array}{rl}
-\Delta u &= \lambda \left(u^{2^*-1}+ \displaystyle \chi_{\{u<a\}}u^{-\de} \right), u  > 0~~\text{in} ~~\Omega,  \\
u & = 0 ~\text{  on   }~
 \partial \Omega.
\end{array}  \end{eqnarray*}

 \noindent We study the existence and the global multiplicity of solutions to the above problem.}}
\end{center}
\section{Introduction}
\def\theequation{1.\arabic{equation}}\makeatother
\setcounter{equation}{0}
\noindent Let $\Omega$ be a bounded domain in $\mathbb R^{N}$, $N\geq 3$ with smooth boundary, $a>0, \lambda>0$ and $0<\de<3$.
 Define $2^*:=\displaystyle\frac{2N}{N-2}$. Denote by $\chi_A$ the characteristic function of a set $A$.
Consider the following elliptic problem with singular and
discontinuous nonlinearity:
\begin{eqnarray*}
 (P_\la^a)~~~~ \Biggl\{\begin{array}{rl}
-\Delta u &= \lambda \left(u^{2^*-1}+ \displaystyle \chi_{\{u<a\}}u^{-\de}\right), u  > 0~~~\text{in} ~~\Omega,  \\
u & = 0 ~\text{  on   }~
 \partial \Omega.
\end{array}  \end{eqnarray*}
\begin{definition}\label{sol-def}
We say that $u\in H^1_0(\Om)$ is a weak solution of $(P_\la^a)$ if ess $\inf_{K} u>0$ for
any compact set $K\subset \Om$ and
\begin{equation}\label{weak soln}
\int_{\Om}\nabla u \nabla \varphi = \la \int_{\Om} (\chi_{\{u<a\}} u^{-\delta}+ u^{2^*-1})\varphi
 \,\,\,\,~ \mbox{ holds for all } \va \in C^{\infty}_0(\Om).
\end{equation}
\end{definition}
\begin{remark} Any solution $u$ of $(P_\la^a)$ belongs to $L^{q}_{loc}(\Om)$
for all $q<\infty$(see for instance, Lemma B.3 of \cite{Str}).Thus $u\in W^{2,q}_{loc}(\Om)$ for all $q<\infty$ and hence
$u\in C^{1,\alpha}_{loc}(\Om)$  for all $\alpha\in(0,1)$.
\end{remark}
\noindent The formal energy functional $E_\la^a(u)$ associated with the problem
$(P_\la^a)$ is given by
\begin{eqnarray*}
E_\la^a(u) = \frac{1}{2}\int_{\Om}|\nabla u|^2\,-\la \int_\Om
G(u) -\frac{\la}{2^{*}} \int_\Om |u|^{2^{*}}
\end{eqnarray*}
where
$$G(u)= \begin{cases}
0 & \mbox{ if } u \leq 0 \\
\displaystyle (1-\de)^{-1}u^{1-\de} & \mbox{ if }0< u<\frac{a}{2},\\
\displaystyle (1-\de)^{-1}(a/2)^{1-\de}+ \int_{a/2}^u
\chi_{\{t<a\}}t^{-\de}dt & \mbox{ if } u\geq a/2,
\end{cases}$$
for $0<\de<3,\de\neq 1$
and for $\de=1$ we replace the terms of the form $(1-\de)^{-1} x^{1-\de}$ in the above definition with the term $\log x$.

We show the following multiplicity result for the problem $(P_\la^a)$.
\begin{theorem}
 For any $a>0$, there exists $\Lambda^a>0$ such that
\begin{enumerate}
\item[(i)]  $(P_\la^a)$ has no solution for any $\la>\La^a$,
\item[(ii)]
$(P_\la^a)$ admits at least two solutions  for  any $\lambda\in(0,\Lambda^a)$.
\end{enumerate}
\end{theorem}

\noi The study of such problems with discontinuous nonlinearities
has  increased remarkably in the last few years due to their occurrence in the modeling of various
physical problems like the obstacle problem,
the seepage surface problem and the Elenbass equation (see \cite{CHANG3},\cite{CHANG2}).
The singular nature of the nonlinearity in $(P_\la^a)$ is motivated by the celebrated
work of Crandall, Rabinowitz and Tartar in \cite{CRT}
which is further studied extensively in
\cite{DPJS},\cite{SS3} \cite{DMO},\cite{GR}, \cite{GS}, \cite{HKS} and \cite{HSS1}.\\
\noi In the pioneering work of Ambrosetti-Brezis-Cerami\cite{ABC},
it was shown that a combination of convex and concave nonlinearities
results in multiple positive solutions for the Dirichlet problem with the model nonlinearity
$\la u^q+ u^\al$, $0<q<1<\al\le\frac{N+2}{N-2}$. In \cite{H} and \cite{HSS}, the  authors have proved similar  multiplicity results
when the nonlinearity in $(P_\la^a)$  has no jump discontinuity and the exponent $\delta$
on the singular term satisfies $0<\delta<1$.
This range for $\delta$ was extended to $0<\delta<3$
 in \cite{AJ} and \cite{DPJS} where the critical and singular nonlinear
problem (again without the jump discontinuity) is discussed in  $\mathbb R^2$.

The problem with jump discontinuity  but without the singular term have been studied in
\cite{AB}, \cite{BT}, \cite{HS} and \cite{SS}.

\noindent In all the above mentioned works, the main methods
used are variational techniques and the generalized gradient theory for locally Lipschitz functionals
as developed in \cite{CHANG} and \cite{CLRK}.
But,  due to the discontinuous and singular
nature of the nonlinear term in our problem, the associated functional is neither differentiable
nor locally Lipschitz in $H^1_0(\Omega)$  and hence  both these techniques can not be used directly.
Therefore, in section 2, we first regularize the discontinuity in $(P_\la^a)$
to make the corresponding functional differentiable and then obtain a first
solution for $(P_\la^a)$ as a limit of the solutions of the regularized problem.
Here we give only the outline of the proof for the existence of the first solution, which
is discussed thoroughly in \cite{SS3}, indicating only the requi modifications.
We then prove that this solution is also a local minimum of
the functional $E^a_\la$ associated with $(P_\la^a)$ in $\HH$ topology.
Since $E^a_\la$ is not in general Fr\'echet differentiable in $\HH$, the ``$H^1$ versus
$C^1"$ result of \cite{BN} can not be used. Instead, an appropriate
use of Hopf's Lemma helps to handle the discontinuity.
In section 3, we prove the existence of a second solution by considering the translate of the problem $(P_\la^a)$
by the first solution and then showing the existence of a solution to the translated problem.
The functional $I_\la$ associated with the translated problem turns out to be locally Lipschitz and hence the theory of generalized gradients
can be applied to prove the existence of the second solution. In this section, we employ Ekeland's variational principle and concentration-compactness ideas
to show the existence of the second positive solution.

\section{Existence of a first solution for $(P_\la^a)$}
\def\theequation{2.\arabic{equation}}\makeatother
\setcounter{equation}{0}
\noindent In this section, we obtain a solution of
$(P_\la^a)$ using the regularizing techniques  similar to that in  \cite{SS3}. Nevertheless, we give an outline of the arguments here for completeness.
Define
\begin{equation}\label{Lambda}
 \La^a=\sup \{\la>0:(P_\la^a) \mbox{ has at least one solution}\},
\end{equation}
and
\begin{eqnarray}\label{eigenfunction}
     \phi_\de =\left\{  \begin{array}{lcc}
   e_1 & \;\;\;\;\;\; & 0<\delta<1,\\
   e_1 \left(-\log e_1\right)^{\frac{1}{2}}& \;\;\;\;\;\; &\delta =1,\\
   e_1^{\frac{2}{\delta+1}} & \;\;\;\;\;\; & 1< \delta,
\end{array}
\right.
\end{eqnarray}
where $ e_1$ is the first positive eigenfunction of $-\Delta $ on $\HH$ 
with $\| e_1\|_{L^\infty(\Om)}$ fixed as a number less than 1.
\begin{lemma}\label{nonexistence}
 $0<\La^a<\infty$.
\end{lemma}
\noindent\textbf{Proof:}
Let $(P_\la^a)$ admit a solution $u_\la$. Since the nonlinearity on the right hand side of $(P_\la^a)$
 is superlinear near infinity, there exists a constant $K=K(a)>0$
such that for all $t>0$, we have  $t^{2^*-1}+\chi_{\{t<a\}}t^{-\de}>Kt$.
Let $\la_1$ be the first eigenvalue of $-\De$ on $\HH$
with the corresponding eigenfunction $ e_1$.
Then multiplying $(P_\la^a)$ by $ e_1$
we get
\begin{align*}
 \la_1\int_\Om u_\la  e_1&=\la \int_\Om (u_\la ^{2^*-1}+\chi_{\{u_\la<a\}}u_\la^{-\de}) e_1 \\
&\geq\la K\int_\Om u_\la  e_1.
\end{align*}
This implies $\La^a<\infty$.
Now we show that $0<\La^a$.
Consider the following  singular problem without the jump-discontinuous term:
\begin{eqnarray*}
 (P_\lambda^{\infty})~~~~ \;\;\;\;\Biggl\{\begin{array}{rl}
-\Delta u &= \lambda ( u^{-\de}+ u^{2^*-1})~~\text{in} ~~\Omega, \\
u & > 0 ~\text{  in  }~ \Omega,\\
u & = 0 ~\text{  on   }~
 \partial \Omega.
\end{array}  \end{eqnarray*}
The existence and multiplicity of solutions of a problem in $\mathbb{R}^2$ analogous to $(P_\la^{\infty})$ has
been studied in \cite{AJ}, \cite{DPJS} and \cite{SS3}.
For $N \geq 3$, a similar approach works as we show now.
From Theorems 1.1,  2.2 and 2.5 in \cite{CRT} 
we can find a unique $v_\la \in \HH$ solving the following purely singular problem for all $\la>0$:
\begin{eqnarray}\label{reg-sing}
\left.\begin{array}{rl}
-\Delta v &= \lambda  v^{-\de},\;\; v>0 ~~\text{in} ~~\Omega, \\
v& = 0 ~\text{  on   }~
 \partial \Omega.
\end{array}  \right\}\end{eqnarray}
It can also be shown (see \cite{CRT} again) that $v_\la\to 0$ uniformly in $\Om$ as $\la \to 0^+$.
Clearly, $v_\la$ is a subsolution to $(P_\lambda^{\infty})$.
 Let $z_\la \in\HH$ solve
\begin{eqnarray}
 \left.\begin{array}{rl}
-\Delta z_\la &=\la,\;\;\; z_\la>0 ~~\text{in} ~~\Omega, \\
z_\la & = 0 ~\text{  on   }~
 \partial \Omega.
\end{array}\right\}  \end{eqnarray}
Define ${\tilde{w}_\la=v_\la+ z_\la}.$ Note that, if $\la_0>0$ is small enough, $\tilde{w}_\la$ is a supersolution to
$(P_\la^\infty)$ for all $\la<\la_0$. Furthermore, $\tilde{w}_\la \to 0$ uniformly on $\Omega$ as $\la \to 0^+$.
 Let $\mathcal{M}_\la= \{u \in \HH : v_\la \leq u \leq \tilde{w}_\la \text{ in } \Om \}$. It is easy to see that $\mathcal{M}_\la$
is a closed convex (hence weakly closed) set in $\HH$.
Now, define the following iterative scheme for all $\la<\la_0$ :
\begin{eqnarray*}
 \;\;\;\;\Biggl\{\begin{array}{rl}
 u_0&=v_\la;\\
-\Delta u_n -\la u_n^{-\de}&= \la u_{n-1}^{2^{*}-1},~
u_n  > 0, ~\text{  in  }~ \Omega,\\
u_n & = 0 ~\text{  on   }~
 \partial \Omega, \;\; n= 1,2,3, \cdots
\end{array}  \end{eqnarray*}
The above scheme is well defined as we can solve for $u_n$ in the closed convex set $\mathcal{M}_\la$
using the Perron's method in variational guise (see \cite{Str}).
As $-\Delta u -\la u^{-\de}$ is a monotone operator in
$\mathcal{M}_\la$, we get that
 $\{u_n\}$ is a non-decreasing sequence.
Thus by standard compactness, we can find $u_\la\in H^1_0(\Om)\cap C^\alpha(\bar\Om)$ for some $\al\in(0,1)$
such that $u_n\rightharpoonup u_\la$ in $H^1_0(\Om)$ and $u_n\rar u_\la$ in $C^\alpha(\bar\Om)$. 
Clearly, since the iteration above started from $v_\la$, we obtain that the solution $u_\la$ obtained is infact a minimal solution.
We note that  $v_\la \leq u_\la \leq \tilde{w}_\la$. 
Also, since $\|u_\la\|_{L^\infty(\Om)}\rar 0$ for $\la \rar 0$, $u_\la$ solves $(P_\la^a)$ for $\la>0$ small and hence
$\La^a>0$.\hfill{\qed}

\begin{corollary}
 From the proof above, it follows that $(P_\la^a)$ admits a solution for $\la>0$ small.
\end{corollary}

\begin{proposition} \label{hark}
Let $\delta>0$.  Given $t_0>0$, there exists a function $p \in C^2((0,t_0)) \cap C([0,t_0])$ satisfying:\\
(i)$ -p^{\prime\prime}=p^{-\delta} $  in $(0,t_0)$,\\
(ii) $p(0)=0$,\\
(iii) $p(s)>0$ in $(0,t_0)$,\\
(iv) For $t$ small, ${p(t) \sim t}$ if $0<\delta < 1$,  $p(t) \geq c t$ if $\delta =1$and $p(t) \sim t^{\frac{2}{\delta+1}}$ if $\delta>1$.
\end{proposition}
\proof Let $0<\delta \leq 1$. Existence of  a  positive function $p \in  C^2((0,t_0+1)) \cap C([0,t_0+1])$ satisfying the equation $ -p^{\prime\prime}=p^{-\delta} $  in $(0,t_0+1)$ and the boundary condition $p(0)=p(t_0+1)=0$  follows from theorem 1.1 in \cite{CRT}. Denote by $\rho_1$ the first (positive) eigenfunction of the interval $(0,t_0+1)$. An easy comparison shows that  $p(t) \geq c \rho_1(t)$ and hence $p(t) \geq ct$ for all small $t>0$ and some $c>0$. If $\delta<1$ it follows that $p^{-\delta} \in L^p((0,t_0+1))$ for some $p>1$ and hence by regularity $p \in C^1([0,t_0+1])$. Therefore, $p(t) \sim t$ near $t=0$ for $\delta<1$.

If $\delta>1$, we take ${p(t)=[\frac{(1+\de)^2}{2(\de-1)}]^{\frac{1}{1+\de}} t^{\frac{2}{\de+1}}}$ for $t>0$.
\hfill\qed

\noindent

We now consider the following two purely singular discontinuous problems:
\begin{eqnarray*}
 (S^a_\lambda)~~~~ \Biggl\{\begin{array}{rl}
-\Delta w &= \la \chi_{\{w<a\}} w^{-\delta},~~w>0~~\text{in} ~~\Omega,  \\
w & = 0 ~\text{  on   }~
 \partial \Omega.
\end{array}  \end{eqnarray*}
and
\begin{eqnarray*}
 (S_\lambda^{a, \epsilon})~~~~ \Biggl\{\begin{array}{rl}
-\Delta w_\e &= \la \chi_\epsilon(w_\e-a)w_\e^{-\delta},~~w_\e>0~~~\text{in} ~~\Omega,  \\
w_\e & = 0 ~\text{  on   }~
 \partial \Omega,
\end{array}  \end{eqnarray*}
where $$\chi_\e(t)=\begin{cases}
                1 &\mbox{ if } t\leq -\e,\\
-(\frac{t}{\e}) & \mbox{ if }-\e<t<0,\\
0& \mbox{ if }t\geq 0.
\end{cases}$$
The existence of $w_\e\in\HH$ solving $(S_\la^{a,\e})$ and satisfying $w_\e\geq c\phi_\de$ ($c$ independent of $\e$) 
follows from proposition \ref{hark} and  theorem 2.2 of \cite{CRT}.
For all $\la>0$, a solution $w_\la$ to $(S^a_\lambda)$
is obtained as the weak limit of the sequence of solutions $\{w_\e\}\subset \HH$ of $(S_\la^{a,\e})$ (for details, see lemma 2.3 in \cite{SS3}).
\begin{theorem}\label{for all lambda}
 $(P_\la^a)$ admits a solution $u_\la$  for all $\la\in(0,\La^a)$. 
 Furthermore, $u_\la$ is global minimum of $E_\la^a$ in the convex set 
 $\mathcal{\overline M}   := \{u \in \HH : w_\la \leq u \leq \overline u \} \subset H^1_0(\Omega)$  where $w_\la$ is a solution to $(S_\la^a)$ and $\overline u$ is a suitable super solution of $(P_\la^a)$.
\end{theorem}
\noindent
\proof We note that $w_\e$ is a subsolution of the
following problem associated to $(P_\la^a)$.
\begin{eqnarray*}
 (P_\lambda^{a,\e})~~~~\;\;\; \Biggl\{\begin{array}{rl}
-\Delta u &= \la  (\chi_\epsilon(u-a)u^{-\delta}+ u^{2^{*}-1}), u>0 ~~\text{in} ~~\Omega,  \\
u & = 0 ~\text{  on   }~
 \partial \Omega.
\end{array}
 \end{eqnarray*}
As before, for  $0< \de < 3,\; \de \neq 1,$ define the following primitive:
$$G_\e(u)= \begin{cases}
0 & \mbox{ if } u \leq 0,\\
\displaystyle (1-\de)^{-1}u^{1-\de} & \mbox{ if }0< u<\frac{a}{2},\\
\displaystyle (1-\de)^{-1}(a/2)^{1-\de}+ \int_{a/2}^u
\chi_\e(t-a)t^{-\de}dt & \mbox{ if } u\geq a/2.
\end{cases}$$
If $\de=1$ we replace the terms of the form $(1-\de)^{-1} x^{1-\de}$ in the above definition with the term $\log x$.
Then the formal energy functional on $\HH$ associated with the problem $(P_\la ^{a,\e})$ is
\begin{eqnarray*}
 E_\la^{a,\e}(u) = \frac{1}{2}\int_{\Om}|\nabla u|^2-\la \int_\Om G_\e(u) -\frac{\la}{2^{*}} \int_\Om |u|^{2^{*}}.
\end{eqnarray*}
Given any  $\la\in(0,\La^a)$, there exists $\bar\la>\la$ such that $(P^a_{\bar\la})$ admits a solution $\overline u$ and
by the definition of $\chi_\e$, $\overline u$ is a supersolution of $(P_\la^{a,\e})$.
Since $-\Delta(w_\e-\overline{u}) \leq \la (\chi_\e(w_\e-a)w_\e^{-\de}
-\chi_\e(\overline{u}-a)\overline{u}^{-\de})$, from the non-increasing nature of
the map $t \mapsto \chi_\e(t-a) t^{-\de},\; t>0,$ we get
 $ w_\e \leq \overline{u}$. Then
the existence of a solution $u_\e$ of $(P_\la^{a,\e})$ is obtained as
local minimizer of $E_\la^{a,\e}$ over the convex set
$\mathcal{M_\e}= \{u\in \HH: w_\e\leq u \leq \overline{u}\}$.
Also, using the same arguments as in \cite{SS3} and \cite{H}, it can be proved that $u_\e$
is a local minimizer of $E_\la^{a,\e}$ in $\HH$.
As $u_\e$ solves $(P_\la^{a,\e})$,
it is easy to check that $\{u_\e\}$ is bounded in $\HH$ and hence weakly converges to some $u_\la \in \HH$.
Then by following the convergence arguments of Lemma 2.3 in \cite{SS3},
it is easy to check that $u_\la$ satisfies $(P_\la^a)$.

 Then as in Proposition 3.1 and Lemma 3.4 of \cite{SS3}, we infer that
 $u_\la$ is a  minimizer of $E_\la^a$ in
 $\mathcal{\overline M} := \{u \in \HH : w_\la \leq u \leq \overline u \}$.
 \hfill{\qed}

\noindent Now we claim that $u_\la$ is a local minimum of $E_\la^a$ in $\HH$.
Here we will follow the same approach as in \cite{SS3} and \cite{SS} and thus be sketchy in the proof.

For $A\subset\mathbb R^N$ we denote $d(x,A)=$ dist$(x,A)$ and by $|A|$ the $N-$dimensional Lebesgue measure of $A$.

\begin{theorem}\label{h1c1}
Let $a>0$. For $\lambda\in(0,\Lambda^a)$, $u_{\lambda}$
is a local minimum for $E_{\lambda}^a$ in $H^{1}_{0}(\Omega)$.
\end{theorem}
\proof
We assume that $u_\la$
is not a local minimum of $E_\la^a$ in $H^1_0(\Om)$ and derive a contradiction.
Let $\{u_n\}\subset\HH$ be such that $u_n\rightarrow
u_\la$ in $\HH$ and $E_\la^a(u_n)<E_\la^a(u_\la)$.
For $\underline u=w_\la$ and solution $\bar u$ of $(P_{\bar\la}^a)$ where $0<\la<\bar\la<\La^a$,
define
$v_n=\max\{\underline{u},\min\{u_n,\bar u\}\},\;
\bar{w}_n=(u_n-\bar u)^+$, $\underline{w}_n=(u_n-\underline{u})^{-}$, $\overline {S}_n=\text{ support }(\bar{w}_n)$
and $\underline{S}_n=\text{ support }(\underline{w}_n).$\\
 \textbf{Claim:} $|\bar S_n|,|\underline S_n|$ and $ \|\bar w_n\|_{\HH}\rightarrow 0$ as $n\rightarrow\infty$.\\
Proof of claim: First to estimate $|\bar S_n|$, we set $\Om_\sigma=\{x\in \Omega:d(x,\partial\Omega)>\sigma\}$ and
$\Om_{\sigma_1}=\{x\in \Omega_\sigma:d(x,\partial\Omega_\sigma)>\sigma_1\}$.
For a given $\e>0$, we choose $\sigma,\sigma_1>0$ sufficiently small
such that $|\Omega \setminus \Omega_\sigma|<\frac{\epsilon}{3}\text{ and }
 |\Omega_\sigma \setminus \Omega_{\sigma_1}|<\frac{\epsilon}{3}.$
First we prove that 
there exists a constant $C>0$ such that
\begin{equation}\label{revised1}
\bar u(x)-u_\la(x)>Cd(x,\partial\Om_\sigma)\text{ and }\underline u+Cd(x,\partial\Om_\sigma)<u_\la \text{ in }\Om_{\sigma_1}.
\end{equation}
\noindent For proving \eqref{revised1}  note that as $\bar u$ is not a solution of $(P_\la^a)$, we have $\overline{u} \not \equiv u_\lambda$ and hence we can choose a small enough  ball $B\subset\subset\Om_{\sigma_1}$
 such that $\bar u\geq u_\la+2\gamma$ in $B$ for some $\gamma>0$.
 Now consider a solution $v$ of the following problem.
 \begin{eqnarray}\label{existence of v}\left. \begin{array}{rl}
  -\Delta v&=\la u_\la^{-\delta}\Psi(v-(\bar u-u_\la))\text{ in }\Omega_{\sigma}\backslash B,  \\
  v&=\gamma \text{ on }\partial B,~v=0\text{ on }\partial\Om_{\sigma}
   \end{array}\right\}
 \end{eqnarray}
 where $\Psi(s)=1$ if $s\leq 0$,$\Psi(s)=-1$ if $s>0$.
Then by the elliptic regularity $v\in W^{2,p}(\Om_{\sigma}\backslash B)\cap
  C^{1,\beta}(\overline{ \Om_{\sigma}\backslash B})$ for some $\beta\in(0,1)$ and for all $p\geq 1$.
 Also taking $v^{-}$ as the test function in the above problem and
 noting that $\bar u\geq u_\la$ in $\Om_{\sigma}\backslash B$, we have $v\geq 0$ in $\Om_{\sigma}\backslash B$.
Furthermore,
 \[-\Delta(\bar u -u_\lambda)\geq\lambda (\bar u^{-\delta} \chi_{\{\bar u<a\}}- u_{\lambda}^{-\delta}\chi_{\{u_\la<a\}})
\geq -\la u_\la^{-\delta} \text{ in }\Omega \]
and $\bar u-u_\la\geq 2\gamma\text{ on }\partial B,\;
 \bar u-u_\la\geq0\text{ on }\partial\Om_{\sigma}$.
Thus,
\begin{equation}\label{new eq 1}
\left.\begin{array}{rl}
 -\Delta(\bar u -u_\lambda-v)&\geq\lambda (- u_{\lambda}^{-\delta} - u_\la^{-\delta} \Psi(v-(\bar u-u_\la))
\text{ in }\Omega_{\sigma}\backslash B,  \\
\bar u-u_\la-v&\geq \gamma\text{ on }\partial B,\\
\bar u-u_\la-v&\geq0\text{ on }\partial\Om_{\sigma}.
\end{array}\right\}
 \end{equation}
Taking $(\bar u-u_\la-v)^{-}$ as the test function in (\ref{new eq 1}) and integrating over $\Om_{\sigma}\backslash B$, we have
\[-\int_{\Om_{\sigma}\backslash B}|\nabla(\bar u-u_\la-v)^-|^2
\geq\la\int_{\Om_{\sigma}\backslash B}[- u_{\lambda}^{-\delta} - u_\la^{-\delta} \Psi(v-(\bar u-u_\la))](\bar u-u_\la-v)^{-}.\]
Now $\Psi(v-(\bar u-u_\la))=-1$ if $(\bar u-u_\la-v)^{-}>0$ and thus the right hand side in the above inequality is zero.
This implies $(\bar u-u_\la-v)^-\equiv 0$ in $\Om_{\sigma}\backslash B$, i.e., $v\leq\bar u-u_\la$ in $\Om_{\sigma}\backslash B$. Therefore,
$-\Delta v=\la u_\la^{-\delta}$ in $\Om_{\sigma}\backslash B$, $v\in C^{1,\beta}(\Om_{\sigma}\backslash B)$, $v>0$
 in $\Om_{\sigma}\backslash B$
and $\frac{\partial v}{\partial\nu}<0$ on $\partial\Om_{\sigma}$ where $\nu$ is the outward unit normal on $\partial\Om_{\sigma}$.
Thus we can find $C>0$ small enough such that $v(x)\geq C d(x,\partial\Om_{\sigma})$ for all $x\in\Om_{\sigma}\backslash B$ 
and hence $(\bar u-u_\la)(x)\geq C d(x,\partial\Om_{\sigma})$ for all $x\in\Om_{\sigma}$. A similar argument can be used to show that $\underline{u}+Cd(x,\partial\Om_\sigma)<u_\la \text{ in }\Om_{\sigma_1}$.
This proves \eqref{revised1}.\\
\noindent Now using \eqref{revised1} we estimate $|\bar S_n|$ as
\begin{align*}
 |\bar S_n|&\leq|\Omega \setminus \Omega_\sigma|+
 |\Omega_\sigma \setminus \Omega_{\sigma_1}|+|\bar S_n\cap\Om_{\sigma_{1}}|\\
&<\frac{\e}{3}+\frac{\e}{3}+\frac{1}{(C\sigma_1)^2}\int_{\bar S_n\cap\Om_{\sigma_{1}}}(u_n-u_\la)^2\\
&\leq\frac{\e}{3}+\frac{\e}{3}+\frac{1}{(C\sigma_1)^2}\|u_n-u_\la\|_{\HH}^2.
\end{align*}
Therefore we get $|\bar S_n|\rar 0$ as $n \to \infty$ and
\begin{align*}
\|\bar w_n\|_{\HH}^2 &=\int_{\bar S_n}|\nabla(u_n-\bar u)|^2
\leq 2\left(\int_{\bar S_n}|\nabla (u_n-u_\lambda)|^2+\int_{\bar S_n}|\nabla(u_\lambda-\bar u)|^2\right)\rightarrow 0
\end{align*}
as $n\rightarrow \infty$.
Using the same approach as above,
we get $|\underline S_n|$ as $n\rightarrow\infty$.
This proves the claim.\\
Note that $v_n\in \mathcal{\overline M}=\{u\in\HH:\underline u\leq u\leq\bar u\}$ and  $u_n=v_n-\underline{w}_n+\bar {w}_n$. Also,
 $E_{\lambda}^a(u_{n})=E_{\lambda}^a(v_{n})+A_{n}+B_{n}$ where
\[A_{n}=\frac{1}{2}\int_{\bar {S}_{n}}(|\nabla u_n|^{2}-|\nabla \bar u|^{2})-
\lambda\int_{\bar {S}_{n}}\left(G(u_n)-G(\bar u)\right)
 -\frac{\la}{2^*}\int_{\bar {S}_{n}}\left(|u_n|^{2^*}-\bar u^{2^*}\right),\]
\[B_{n}=\frac{1}{2}\int_{\underline S_{n}} (|\nabla u_n|^{2}-|\nabla \underline u|^{2})-
\lambda\int_{\underline {S}_{n}} \left(G(u_n)-G(\underline u)\right)
 -\frac{\la}{2^*}\int_{\underline {S}_{n}}\left(|u_n|^{2^*}-\underline u^{2^*}\right).\]
As $u_\la$ is minimizer of $E_\la^a$ over $\mathcal{\overline M}$ (see theorem 2.1) and
$v_n\in \mathcal{\overline M}$ we have $E_\la^a(u_n) \geq
E_\la^a(u_\la)+A_n+B_n$. Now we claim that $A_n, B_n \geq 0$ for all large $n$
which is a contradiction to our assumption that  $E_\la^a(u_n)<E_\la^a(u_\la)$ for all $n$.
Note that
\begin{align*}
 A_n&=\frac{1}{2}\int_{\bar {S}_{n}}(|\nabla u_n|^{2}-|\nabla \bar u|^{2})-
\lambda\int_{\bar{S}_n}\left(G(u_n)-G(\bar u)\right)
 -\frac{\lambda}{2^*}\int_{\bar {S}_{n}}\left(|u_n|^{2^*}-\bar u^{2^*}\right)\\
&=\frac{1}{2}\int_{\bar {S}_{n}}|\nabla \bar w_n|^2+\int_{\bar {S}_{n}}\nabla \bar u\cdot\nabla \bar w_n
-\lambda\int_{\bar {S}_{n}}\left(G(\bar u+\bar w_n)-G(\bar u)\right)\\
&~~~ -\frac{\lambda}{2^*}\int_{\bar {S}_{n}}\left((\bar u+\bar w_n)^{2^*}-\bar u^{2^*}\right)\\
&\geq\frac{1}{2}\int_{\bar {S}_{n}}|\nabla \bar w_n|^2
+\lambda\int_{\bar {S}_{n}}\left(\chi_{\{\bar u<a\}}\bar u^{-\delta}\bar w_n-(G(\bar u+\bar w_n)-G(\bar u))\right)\\
&~~~+\la\int_{\bar {S}_{n}}\left(\bar u^{2^*-1}\bar w_n-\frac{1}{2^*}((\bar u+\bar w_n)^{2^*}-\bar u^{2^*})\right).
\end{align*}
Now by dividing $\bar S_n$ into three subdomains, viz.,
$\bar S_n\cap\{x\in\Om:a<\bar u(x)\}$,
$\bar S_n\cap\{x\in\Om:\bar u(x)\leq a\leq (\bar u+\bar w_n)(x)\}$ and $\bar S_n\cap\{x\in\Om:(\bar u+\bar w_n)(x)<a\}$,
one can check that the second integral in the
right hand side of the above inequality is nonnegative. Also by the mean value theorem,
for some $\theta=\theta(x)\in(0,1)$ and appropriate positive constants $c_1,c_2,c_3$ we have
\begin{align}\label{A2}
 \int_{\bar {S}_{n}}\bar u^{2^*-1}\bar w_n
-\frac{1}{2^*}\left((\bar u+\bar w_n)^{2^*}-\bar u^{2^*}\right)&=
-\int_{\bar {S}_{n}}((\bar u+\theta\bar w_n)^{2^*-1}-\bar u^{2^*-1})\bar w_n\nonumber\\
&\geq -\int_{\bar {S}_{n}}(\bar u+\bar w_n)^{2^*-2}\bar w_n^2\nonumber\\
&\geq-c_1\int_{\bar {S}_{n}}(\bar u^{2^*-2}+\bar w_n^{2^*-2})\bar w_n^2\nonumber\\
 &\geq-c_2\left(\int_{\bar {S}_{n}}\bar u^{2^*}\right)^{\frac{2^*-2}{2^*}}\|\bar w_n\|_{\HH}^2\nonumber\\
 &~~~~-c_3\|\bar w_n\|_{\HH}^{2^*}.
\end{align}
Thus using (\ref{A2}) we have the following estimation for $A_n$:
\begin{equation}\label{A3}
A_n\geq\frac{1}{2}\|\bar w_n\|_{\HH}^2-\la c_2\;\left(\int_{\bar {S}_{n}}\bar u^{2^*}\right)^{\frac{2^*-2}{2^*}}\|\bar w_n\|_{\HH}^2
-\la c_3\|\bar w_n\|_{\HH}^{2^*}.\end{equation}
Also following the arguments as in \cite{SS3} and using (\ref{A2}) we estimate $B_n$ as
\begin{equation}\label{A4}
B_n\geq\frac{1}{2}\|\underline w_n\|_{\HH}^2-C\left(\int_{\underline{S}_n}\underline u^{2^*}\right)^{\frac{2^*-2}{2^*}}\|\underline w_n\|_{\HH}^2.
\end{equation}
Since $|\bar S_n|,|\underline S_n| $ and $ \|\bar w_n\|_{\HH}\rightarrow 0$ as $n\rightarrow\infty$, we get $A_n,B_n\geq0$.
This completes the proof of the theorem.\hfill{\qed}
\section{Existence of the second solution for $(P_{\la}^a)$}
\def\theequation{3.\arabic{equation}}\makeatother
\setcounter{equation}{0}
\noindent In this section we obtain a second solution for $(P_\la^a)$ for $\la \in (0,\La^a)$
by translating the problem to the solution $u_\la$ obtained in the previous section.
We consider the following translated problem $(\tilde{P}_\lambda^a):$
\begin{eqnarray*}
 (\tilde{P}_\lambda^a)~~~~ \left\{\begin{array}{rl}
-\Delta u &= \lambda\displaystyle \left(\chi_{\{u+u_\la<a\}}(u+u_\la)^{-\de}-\chi_{\{u_\la<a\}}u_\la^{-\de}\right)\\
 &~~~~~~+\la \left(  (u+u_\la)^{2^*-1}-u_\la^{2^*-1}\right)~~\text{in} ~~\Omega, \\
u&>  0\;\text{ in } \Om,\\
u& = 0 ~\text{  on   }~
 \partial \Omega.
\end{array} \right. \end{eqnarray*}
\begin{remark}\label{rem 2 sol}
It is easy to see that if $v_\la \in \HH$ weakly solves $(\tilde{P}_\lambda^a)$,
then $u_\la+v_\la$ weakly solves $(P_\la^a)$.
\end{remark}
\noindent Let us define, for $x\in\Om$,
$$\begin{array}{lll}
 \tilde g(x,s) &= &\displaystyle{( \chi_{\{s+u_\la(x)<a\}}(s+u_\la(x))^{-\de}- \chi_{\{u_\la(x)<a\}}u_\la(x)^{-\de}
 ) \chi_{\mathbb{R}^+}(s)} \, \mbox{\;\;and}\\[2 mm]
 \tilde f(x,s)&=&\displaystyle{\left((s+u_\la(x))^{2^*-1}-(u_\la(x))^{2^*-1}\right)\chi_{\mathbb{R}^{+}}(s).}
\end{array}$$
Let $\tilde{G}(x,t)=\int_0^{t} \tilde g(x,s) ds \,$ and
$\tilde{F}(x,t)=\int_0^{t} \tilde f(x,s) ds.$
Let $\II : \HH \to \RR$ be the energy functional associated with $(\tilde{P}_\la^a)$ defined as below:
\begin{equation}\label{2sol1}
 I_\la(u) = \frac{1}{2}\int_{\Om}|\nabla u|^2-\la \int_\Om \tilde{G}(x,u(x))dx -\la \int_\Om \tilde{F}(x,u(x)) dx.
\end{equation}
\begin{proposition}\label{I Lip}
$I_\la$ is  locally Lipschitz on $\HH$.
\end{proposition}
\noindent
\proof Note that as $\HH \ni u \mapsto \frac{1}{2}\int_\Om|\nabla u|^2-\la \int_\Om \tilde{F}(x,u(x))dx$
is a $C^1$ map, it is sufficient to prove that the
map $:\HH \ni u \mapsto  \int_\Om \tilde{G}(x,u(x))dx
\in \RR$ is  Lipschitz.
 We have,
$$\displaystyle \left| \int_{\Om} \tilde{G}(x,u+v)-\tilde{G}(x,u) \,  \right|
\leq \displaystyle\int_{\Om}\left| \int_{u(x)}^{(u+v)(x)}\tilde{g}(x,s) \,ds \,\right| \,
\leq2\displaystyle\int_{\Om}u_\la^{-\de} |v|.$$
Since $u_\la (x)\geq \underline{u} := w_\la \geq k_1 \phi_\de \geq C d(x,\partial\Om)^{\frac{2}{1+\de}}$ for some $k_1,C>0$ 
where $w_\la$ is a solution to $(S_\la^a)$ (see the proof of Theorem \ref{for all lambda}) and the remarks immediately above this theorem), 
 thanks to Hardy's inequality, it can be easily checked that
$$\int_{\Om}u_\la^{-\de} |v| \, dx  \leq \int_{\Om}\frac{|v|}{Cd(x,\partial\Om)^{\frac{2\de}{1+\de}}}
\, dx=\int_{\Om}\frac{|v|}{C d(x,\partial\Om)} d(x,\partial\Om)^{\frac{1-\de}{1+\de}}\leq
C_1 \, ||v||_{\HH}.
$$
Hence, $I_\la$ is locally Lipschitz.\hfill \qed
\begin{definition}
Let $I : \HH \rar \RR$ be a locally Lipschitz map. The generalized derivative of $I$ at $u $ in the
direction of $\phi$ $($denoted by $I^0(u,\phi))$ is defined as:
$$I^0(u,\phi) =\limsup_{h\rar 0,t \downarrow 0} \frac{I(u+h+t\phi)-I(u+h)}{t}; \;\; u, \phi \in \HH.$$\\
We say that $u$ is a ``generalized" critical point of $I$ if
$I^0(u,\phi) \geq 0$ for all $\phi \in \HH$.
 See \cite{CHANG} (page 103) for more details.
\end{definition}

\begin{remark}\label{w-existence}
From Lemma 4.1 of \cite{SS3}, for $u\geq 0$ and $\phi\in\HH$, we have the following inequality:
\begin{equation}\label{2sol2}
I^0_\la(u,\phi)\leq\int_{\Om} \nabla (u_\la+u) \cdot \nabla \phi \, - \la \int_{\Om}(u_\la+u)^{2^*-1} \phi\,  -\la \int_{\Om}w^\phi(u_\la+u)^{-\de} \phi
\end{equation}
for some measurable function $w^\phi\in[\chi_{\{u_\la+u<a\}},\chi_{\{u_\la+u\leq a\}}]$.
\end{remark}
\begin{remark}\label{measure-zero}
 From Remark 4.4 of \cite{SS3}, suppose for some nontrivial, nonnegative $v_\la \in \HH$
we have $I_\la^0(v_\la, \phi ) \geq 0$ for all $\phi \in \HH$,
 i.e., $v_\la$ is a ``generalized" critical point of $I_\la$. Then, from $(\ref{2sol2})$,
 \begin{equation}\label{eqn2.4}
 \la (u_\la+v_\la)^{2^*-1} \leq -\Delta (u_\la+v_\la) \leq \la [(u_\la+v_\la)^{2^*-1}+ (u_\la+v_\la)^{-\de}]
\end{equation}
in the weak sense.
Let us show \eqref{eqn2.4}. Indeed, as $v_\la\geq 0$ and $I_\la^0(v_\la, \phi)\geq 0$, using (3.2), we have for all $\phi\in H^1_0(\Om),$
 \begin{equation}\label{exp13}
0\leq I^0_\la(v_\la,\phi)\leq\int_{\Om} \nabla (u_\la+v_\la) \cdot \nabla \phi \, - \la \int_{\Om}(u_\la+v_\la)^{2^*-1} \phi\,  -\la \int_{\Om}w^\phi(u_\la+v_\la)^{-\de} \phi.
\end{equation}
Let $\phi\geq 0,$ then we have from (\ref{exp13}), 
$$ \int_{\Om} \nabla (u_\la+v_\la) \cdot \nabla \phi\geq \la \int_{\Om}(u_\la+v_\la)^{2^*-1} \phi\,  
+\la \int_{\Om}w^\phi(u_\la+v_\la)^{-\de} \phi. $$
Since $w^\phi\geq 0$ and given that $\phi\geq 0,$ we have 
\begin{equation} \label{exp13.1}
\int_{\Om} \nabla (u_\la+v_\la) \cdot  \nabla \phi\geq \la \int_{\Om}(u_\la+v_\la)^{2^*-1} \phi 
\end{equation}
or in other words, $-\Delta (u_\la+v_\la)\geq \la (u_\la+v_\la)^{2^*-1}$ in the weak sense. 

Next let us consider a $\phi\in H^1_0(\Om)$ which is non-positive, so that $\psi=-\phi\geq 0.$ 
Again using  $(\ref{exp13})$ we have,
$$ \int_{\Om} \nabla (u_\la+v_\la) \cdot \nabla(- \psi)\geq \la \int_{\Om}(u_\la+v_\la)^{2^*-1}(- \psi)\,  
+\la \int_{\Om}w^{-\psi}(u_\la+v_\la)^{-\de}(- \psi). $$
Multiplying by $-1$ on both sides and using the fact that $w^{-\psi}\in [0,1]$ we get, 
$$ \int_{\Om} \nabla (u_\la+v_\la) \cdot \nabla \psi\leq \la \int_{\Om}(u_\la+v_\la)^{2^*-1} \psi\,  
+\la \int_{\Om}(u_\la+v_\la)^{-\de} \psi . $$
Since $\psi=-\phi$ is any arbitrary non-negative function in $H^1_0(\Om),$ the previous expression implies 
\begin{equation}\label{exp13.2} -\Delta (u_\la+v_\la) \leq \la (u_\la+v_\la)^{2^*-1}+ \la (u_\la+v_\la)^{-\delta} \mbox{ in the weak sense.}\end{equation}
From (\ref{exp13.1}) and (\ref{exp13.2}) we conclude that 
\begin{equation}\label{exp13.4}
 \la (u_\la+v_\la)^{2^*-1} \leq -\Delta (u_\la+v_\la) \leq \la [(u_\la+v_\la)^{2^*-1}+ (u_\la+v_\la)^{-\de}].
\end{equation}
{Note that $-\Delta(u_\la+v_\la)$ is a positive distribution and hence it is given by a positive, regular
  Radon measure say $\mu$. Then using (3.7) we can show that $\mu$ is absolutely continuous with respect to the Lebesgue measure.  Now by Radon Nikodyn theorem there 
  exists a locally integrable function $g$ such that $-\Delta(u_\la+v_\la)=g$ and hence $g\in L^p_{loc}(\Om)$ for some $p>1.$ }
Now using Lemma B.3 of \cite{Str} and elliptic regularity we can conclude that $u_\la+v_\la\in W^{2,q}_{loc}(\Om)$ 
for all $q<\infty$
and for almost every $x \in \Om,$ $-\Delta (u_\la+w_\la)(x)= g(x)\geq \la(u_\la+v_\la)^{2^*-1}(x)>0.$ In particular,
\begin{equation}\label{exp13.5}
-\Delta (u_\la+v_\la)>0 \mbox{ for a.e on } \{ x\in \Om: (u_\la+v_\la)(x)=a\}.   
\end{equation}
On the other hand, 
we have
$-\Delta(u_\la+v_\la)=0$ a.e on the set $\{x: (u_\la+v_\la)(x) = a\}.$ 
This contradicts (\ref{exp13.5}) unless the Lebesgue measure of the set $\{x: (u_\la+v_\la)(x) = a\}$ is zero. 

Therefore, $w^{\phi}= \chi_{\{u_\la+v_\la<a\}}$ a.e. in $\Om$ for any $\phi \in \HH$ and
 hence $u_\la+v_\la$ is a second solution for $(P_\la^a)$.
\end{remark}
\begin{remark}\label{I and E}
 Note that as $I_\la(u)=E^a_\la(u^++u_\la)-E^a_\la(u_\la)+\frac{1}{2}\int_{\Om} |\nabla u^-|^2$
 for any $u \in \HH$ and $u_\la$ is a local minimum of $E^a_\la$ in $\HH$, it follows that $0$ is a local minimum
of $I_\la$ in $\HH$.
 \end{remark}
\noindent Using the Mountain Pass theorem and Ekeland variational principle
we show the existence of a generalized critical point for $I_\la$ which yields a second solution to $(P_\la^a)$.
The method of the proof is along lines similar to those of \cite{SS3}.
Let us define $H^+=\{u\in \HH: u \geq 0 \, \mbox{ a.e in }\Om\}$. Since 0 is a local
minimum of $\II$, there exists a $\rho_0>0$ such that $\II(0) \leq \II(u)$ for $||u||_{\HH}\leq \rho_0$.
The following two cases arise:\\
\textbf{1.}$ZA$ (Zero altitude): $\inf\{ \II(u): ||u||_{\HH}=\rho, u\in H^+\}=\II(0)=0$ for all $\rho \in (0,\rho_0).$\\
\textbf{2.}$MP$ (Mountain Pass): There exists $\rho_1\in (0,\rho_0)$ such that $\inf\{ \II(u): ||u||_{\HH}=\rho_1, u\in H^+\}>\II(0).$
\begin{lemma}\label{ZA}
Let $ZA$ hold for some $\la \in (0, \La^a)$.
Then there exists a nontrivial ``generalized" critical point $v_\la \in H^+$ for $I_\la$.
\end{lemma}
\noindent
 \proof Fix $\rho \in (0,\rho_0)$. Then
there exists a sequence $\{z_n\}\subset H^+$ with $\|z_n\|_{\HH}=\rho$ and $\II(z_n) \leq 1/n.$ Fix $0<r<
\frac{1}{2}\min\{\rho_0-\rho,\rho\}$ and define
$R=\{u\in H^+: \rho-r \leq \|u\|_{\HH}\leq \rho+r\}.$ Clearly $R$ is closed and $\II$ is Lipschitz continuous on $R$ from Proposition \ref{I Lip}.
Thus by Ekeland's variational
principle there exists $\{v_n\}\subset R$ such that
\begin{enumerate}
\item[(i)] $I_{\lambda}(v_{n})\leq I_{\lambda}(z_{n})\leq\displaystyle \frac{1}{n}$,
\vspace{-0.2cm}
\item[(ii)] $\|z_n-v_n\|_{\HH}\leq\displaystyle\frac{1}{n}$ \;\;and
\vspace{-0.2cm}
\item[(iii)] $I_{\lambda}(v_{n})\leq I_{\lambda}(v)+\displaystyle\frac{1}{n}\|v-v_n\|_{\HH}$ for all $v\in R$.
\end{enumerate}
 We note that 
\be \label{dov}
\rho - \frac{1}{n} = \|z_n\|_{\HH} - \frac{1}{n} \leq \|v_n\|_{\HH} \leq \|z_n\|_{\HH} + \frac{1}{n} = \rho + \frac{1}{n}.
\ee

Therefore, for $\xi\in H^+$ we can choose $\epsilon>0$ sufficiently small such that $v_n+\epsilon(\xi-v_n)\in R$ for all large $n$.
Then by $(iii)$ we get
$$\displaystyle\frac{I_{\lambda}(v_n+\epsilon(\xi-v_n))-I_{\lambda}(v_n)}{\epsilon}\geq-\frac{1}{n}\|\xi-v_n\|_{\HH}.$$
Letting $\epsilon\rightarrow 0^+$, we conclude
$$I_{\lambda}^{0}(v_n,\xi-v_n)\geq-\displaystyle\frac{1}{n}\|\xi-v_n\|_{\HH}\text{ for all } \xi\in H^+.$$
From Remark \ref{w-existence}, for any $\xi\in H^+$, there exists $w^{\xi-v_n}_n\in[\chi_{\{u_\la+v_n<a\}},\chi_{\{u_\la+v_n\leq a\}}]$ such that
\begin{align}\label{2sol3}
\int_{\Om} \nabla (u_\la+v_n)&\cdot\nabla (\xi-v_n) - \la \int_{\Om}(u_\la+v_n)^{2^*-1}(\xi-v_n)\nonumber\\
&-\la \int_{\Om}w^{\xi-v_n}_n(u_\la+v_n)^{-\de}(\xi-v_n)\geq-\frac{1}{n}\|\xi-v_n\|_{\HH}.
\end{align}
Since $\{v_n\}$ is bounded in $H^{1}_{0}(\Omega)$, we may assume $v_n \rightharpoonup v_\lambda\in H^+$ weakly in $H^{1}_{0}(\Omega)$.
Now by following the same arguments as in Lemma 4.2 of \cite{SS3} we can show that $v_\la$ is a generalized critical point for $I_\la.$
It remains to show that $v_\la \not \equiv 0.$ Note that if $I_\la(v_\la)\neq 0$ we are done. So assume $I_\la(v_\la)=0$.
Since $\|v_n\|_{\HH}\geq \rho/2$ for all large $n$ (see \eqref{dov}), it is sufficient to show that $v_n \rar v_\la$ strongly in
$\HH$. Taking $\xi=v_\la$ in (\ref{2sol3}) we get
\begin{align}\label{2sol10}
\int_\Om&\nabla(u_\la+v_\la)\cdot\nabla(v_\la-v_n)-\la\int_\Om (u_\la+v_n)^{2^*-1}(v_\la-v_n)\nonumber\\
&-\la\int_\Om w_n^{v_\la-v_n}(u_\la+v_n)^{-\de}(v_\la-v_n)+\frac{1}{n}\|v_\la-v_n\|_{\HH}
\geq\|v_\la-v_n\|_{\HH}^2.
 \end{align}
For any measurable set $E \subset \Om$, as $u_\la\geq k_1\phi_\de$ and $v_n\in H^+$, thanks to Hardy's inequality, we have
\begin{align}\label{Hardy}
 \int_{E}  w_{n}^{v_\la-v_n}|v_n-v_\la|(v_n+u_\la)^{-\de}
 &\leq  \displaystyle C \int_{E}\frac{|v_n-v_\la|}{u_\la^{\de}}\nonumber\\
 &\leq \displaystyle C \int_{E}\frac{|v_n-v_\la|}{d(x,\partial\Om)^{\frac{2\de}{1+\de}}}\nonumber\\
&\leq  \displaystyle C \int_{E}\frac{|v_n-v_\la|}{d(x,\partial\Om)} d(x,\partial\Om)^{\frac{1-\de}{1+\de}}\nonumber\\
&\leq C \|v_n-v_\la\|_{\HH}\, \|d(x,\partial\Om)^{\frac{1-\de}{1+\de}}\|_{L^2(E)}.
\end{align}
Since $v_n\rar v_\la$ pointwise a.e. in $\Om$,
by Vitali's convergence theorem,
\begin{equation}\label{2sol11}
 \displaystyle\int_{\Om}  w_{n}^{v_\la-v_n}|v_n-v_\la|(v_n+v_\la)^{-\de} \,  \to 0 \text{ as }n\rar \infty.
\end{equation}
Also from Brezis-Lieb lemma (\cite{BL}) we have
\begin{align}\label{2sol12}
\int_\Om (u_\la+v_n)^{2^*-1}(v_\la-v_n)&=\int_\Om (u_\la+v_n)^{2^*-1}(u_\la+v_\la)-\int_\Om (u_\la+v_n)^{2^*}\nonumber\\
&=-\|v_\la-v_n\|_{L^{2^*}(\Om)}^{2^*}+o_n(1).
\end{align}
Now using (\ref{2sol11}) and (\ref{2sol12}) in (\ref{2sol10}) we get
\begin{equation}\label{new0}
 \|v_\la-v_n\|_{\HH}^2-\la\|v_\la-v_n\|^{2^*}_{L^{2^*}(\Om)}\leq o_n(1).
\end{equation}
Also taking $\xi=2v_n$ in (\ref{2sol3}) and using the fact that $u_\la$ solves $(P_\la^a)$
we get
\begin{eqnarray*}
-\frac{1}{n}\|v_n\|_{\HH}&\leq& \int_\Om\nabla u_\la\cdot\nabla v_n
+\int_\Om|\nabla v_n|^2-\la\int_\Om (u_\la+v_n)^{2^*-1}v_n\nonumber\\
&&-\la \int_\Om w_n^{v_n}v_n(u_\la+v_n)^{-\delta}\nonumber\\
&=&\|v_n\|_{\HH}^2-\la\int_\Om ((u_\la+v_n)^{2^*-1}-u_\la^{2^*-1})v_n\nonumber\\
&&+\la\int_\Om\left(\chi_{\{u_\la<a\}}u_\la^{-\delta}-w_n^{v_n}(u_\la+v_n)^{-\delta}\right)v_n\nonumber\\
&=& \|v_\la\|_{\HH}^2+\|v_n-v_\la\|_{\HH}^2-\la\int_\Om\tilde f(v_n)v_n\nonumber\\
&&+\la\int_\Om\left(\chi_{\{u_\la<a\}}u_\la^{-\delta}-w_n^{v_n}(u_\la+v_n)^{-\delta}\right)v_n+o_n(1). 
\end{eqnarray*}
Now as $v_\la$ solves $(\tilde P_\la^a)$ we have
$$\|v_\la\|_{\HH}^2=\la\int_\Om\tilde f(v_\la)v_\la +
\la\int_\Om \left(\chi_{\{u_\la+v_\la<a\}}(u_\la+v_\la)^{-\delta}-\chi_{\{u_\la<a\}}u_\la^{-\delta}\right)v_\la.
$$
Using this identity in above inequality  we get,
\begin{align}\label{new1}
-\frac{1}{n}\|v_n\|_{\HH}&\leq
\|v_n-v_\la\|_{\HH}^2-\la\int_\Om\left(\tilde f(v_n)v_n-\tilde f(v_\la)v_\la \right)\nonumber\\
&~~~~~~+\la\int_\Om \left(\chi_{\{u_\la+v_\la<a\}}(u_\la+v_\la)^{-\delta}-\chi_{\{u_\la<a\}}u_\la^{-\delta}\right)v_\la\nonumber\\
&~~~~~~+\la\int_\Om\left(\chi_{\{u_\la<a\}}u_\la^{-\delta}-w_n^{v_n}(u_\la+v_n)^{-\delta}\right)v_n+o_n(1).
\end{align}
Using again Brezis-Lieb lemma it is easy to check that
\[\int_\Om\tilde f(v_n)v_n-\tilde f(v_\la)v_\la=\|v_n-v_\la\|^{2^*}_{L^{2^*}(\Om)}+o_n(1).\]
Also as $v_n\rar v_\la$ pointwise a.e. in $\Om$ and $|\{x\in\Om:(u_\la+v_\la)(x)=a\}|=0$,
using estimates similar to the one in \eqref{Hardy} we have
\begin{align*}
&\int_\Om \left(\chi_{\{u_\la+v_\la<a\}}(u_\la+v_\la)^{-\delta}-\chi_{\{u_\la<a\}}u_\la^{-\delta}\right)v_\la\\
&~~~~~~~~+\int_\Om\left(\chi_{\{u_\la<a\}}u_\la^{-\delta}-w_n^{v_n}(u_\la+v_n)^{-\delta}\right)v_n
=o_n(1).
\end{align*}
Thus (\ref{new1}) implies
\begin{align}\label{new2}
 o_n(1)&\leq\|v_n-v_\la\|_{\HH}^2-\la\|v_n-v_\la\|^{2^*}_{L^{2^*}(\Om)}.
\end{align}
Also as $I_\la(v_n)\leq\displaystyle\frac{1}{n}$
and $\tilde F(v_n)=\displaystyle\frac{(u_\la+v_n)^{2^*}}{2^*}-\displaystyle\frac{u_\la^{2^*}}{2^*}-u_\la^{2^*-1}v_n$,
we have
\begin{align*}I_\la(v_n)&=\frac{1}{2}\|v_n\|_{\HH}^2-\frac{\la}{2^*}\|u_\la+v_n\|_{L^{2^*}(\Om)}^{2^*}
+\frac{\la}{2^*}\|u_\la\|_{L^{2^*}(\Om)}^{2^*}\\
&\quad\quad+\la\int_\Om u_\la^{2^*-1}v_n
-\la\int_\Om\tilde G(v_n)\\
&\leq\frac{1}{n}.\end{align*}
From the fact that $v_n\rightharpoonup v_\la$ weakly in $\HH$, this implies
\begin{equation}\label{new4}
 \frac{1}{2}\|v_n-v_\la\|_{\HH}^2-\frac{\la}{2^*}\|v_n-v_\la\|_{L^{2^*}(\Om)}^{2^*}+I_\la(v_\la)
+\la\int_\Om\tilde G(v_\la)-\la\int_\Om\tilde G(v_n)
\leq o_n(1).
\end{equation}
Now using the Hardy's inequality and Vitali's convergence theorem as in \eqref{Hardy} one can check that
$\int_\Om\tilde G(v_n)\rar\int_\Om\tilde G(v_\la)$ as $n\rar\infty$. Also as $I_\la(v_\la)= 0$, (\ref{new4}) implies
\begin{equation}\label{new5}
  \frac{1}{2}\|v_n-v_\la\|_{\HH}^2-\frac{\la}{2^*}\|v_n-v_\la\|_{L^{2^*}(\Om)}^{2^*}\leq o_n(1).
\end{equation}
Now from (\ref{new0}), (\ref{new2}) and (\ref{new5}) we get
$\displaystyle(\frac{1}{2}-\displaystyle\frac{1}{2^*})\|v_n-v_\la\|_{\HH}^2\leq o_n(1)$ and hence
$v_n\rightarrow v_\la$ in $H^1_0(\Om)$.\hfill{\qed}\\
\noindent Next we consider the case (MP). As the nonlinearity is critical, we
use the following Talenti functions to study the critical level:
$$V_\e(x)=\frac{C_N\e^{(N-2)/2}}{(\e^2+|x|^2)^{(N-2)/2}},\;\;\;C_N,\e>0.$$
Fix any $y\in \Om_a=\{x\in\Om\;:u_\la(x)<a\}$.
Choose $\eta\in C^\infty_c(\Om)$ such that $0\leq\eta\leq 1$ and $\eta\equiv 1$
on $\overline{B_{r}(y)}$ where $r>0$ is chosen small enough such that
$\overline{B_{r}(y)}\subset \Om_a$. Define
$U_\e(x)=\eta(x)V_\e(x-y).$
Then  as $\e\rar 0$, a standard computation (see \cite{BNcritical}) gives
\begin{equation}\label{critical1}
\int_\Om|U_\e|^{2^*}=\int_{\mathbb R^N}\left\vert V_1\right\vert^{2^*}+o(\e^N)=A+o(\e^{N})
\end{equation}
and
\begin{equation}\label{critical2}
\int_\Om|\nabla U_\e|^2=\int_{\mathbb R^N}|\nabla V_1|^2+o(\e^{N-2})=B+o(\e^{N-2}).
\end{equation}
We have the following lemma.
\begin{lemma}\label{L2}
 There exist $\e_0>0$ and $R_0\geq 1$ such that
\begin{itemize}
 \item [(i)]$I_\la(RU_\e)<I_\la(0)=0$ for all $\e\in(0,\e_0)$ and $R\geq R_0$.
\item[(ii)] $I_\la(tR_0U_\e)<\displaystyle\frac{S^{\frac{N}{2}}}{N\la^{(N-2)/2}}$ for all $t\in(0,1], \e\in(0,\e_0)$
where  $S=\displaystyle\frac{B}{A^{2/2^*}}$ is the best constant of the Sobolev embedding.
\end{itemize}
\end{lemma}
\noindent
\proof Noting that for $v\in H^+$, $E_\la^a(u_\la+v)=E_\la^a(u_\la)+I_\la(v)$, this is equivalent to show that
\begin{itemize}
\item [(i)]$E_\la^a(u_\la+RU_\e)<E_\la^a(u_\la)$ for all $\e\in(0,\e_0)$ and $R\geq R_0$.
\item[(ii)] $E_\la^a(u_\la+tR_0U_\e)<E_\la^a(u_\la)+\displaystyle\frac{S^{\frac{N}{2}}}{N\la^{(N-2)/2}}$ for all $t\in(0,1],\e\in(0,\e_0)$.
\end{itemize}
Now using the fact that $u_\la$ solves $(P_\la^a)$, first we estimate $E_\la^a(u_\la+tRU_\e)$, $t>0$, as follows.
\begin{align*}
E_\la^a(u_\la+tRU_\e)&=\frac{1}{2}\int_\Om|\nabla(u_\la+tRU_\e)|^2-\la\int_\Om G(u_\la+tRU_\e)\\
&\quad -\frac{\la}{2^*}\int_\Om(u_\la+tRU_\e)^{2^*}\\
&=\frac{1}{2}\int_\Om|\nabla u_\la|^2+\frac{R^{2}t^{2}}{2}\int_\Om|\nabla U_\e|^2+tR\int_\Om\nabla u_\la \cdot \nabla U_\e
\\
&\quad -\la\int_\Om G(u_\la+tRU_\e)-\frac{\la}{2^*}\int_\Om(u_\la+tRU_\e)^{2^*}\\
&=\frac{1}{2}\int_\Om|\nabla u_\la|^2+\frac{R^{2}t^{2}}{2}\int_\Om|\nabla U_\e|^2 +\la t R\int_\Om(\chi_{\{u_\la<a\}}u_\la^{-\delta}+ u_\la^{2^*-1})U_\e \\
&~~~-\la\int_\Om G(u_\la+tRU_\e)
-\frac{\la}{2^*}\int_\Om(u_\la+tRU_\e)^{2^*}.
\end{align*}
Now we estimate the critical term $\int_\Om(u_\la+tRU_\e)^{2^*}$  using the one-dimensional inequality in Lemma 4 of \cite{BNminimization} as:
 \begin{align*}
 \int_\Om(u_\la+tRU_\e)^{2^*}&=\int_\Om u_\la^{2^*}+(tR)^{2^*}\int_\Om U_\e^{2^*}+2^*tR\int_\Om u_\la^{2^*-1}U_\e\\
 &~~~~+2^*(tR)^{2^*-1}\int_\Om u_\la U_\e^{2^*-1}+(R_\e+S_\e).
 \end{align*}
 The terms $R_\e$ and $S_\e$ are given by the following expressions:
 \begin{equation}\label{BN0}
 R_\e=\begin{cases}
 O_\e(\int_{\Omega} U_\epsilon^2) &\text{ if } N < 6,\\
 \displaystyle O_\e( \int_{\{u_\la\geq tR U_\e\}}u_\la (tR U_\e)^{2^*-1})&  \text{ if }N \geq 6,
 \end{cases}
 \end{equation}
 and
 \begin{equation}\label{BN01}
 S_\e=\begin{cases}
 O_\e(\int_{\Omega}U_\epsilon^{2^*-2}) &\text{ if } N < 6,\\
 \displaystyle O_\e(\int_{\{u_\la\leq tR U_\e\}}u_\la^{2^*-1} (tR U_\e)) &  \text{ if }N \geq 6.
 \end{cases}
 \end{equation}
  Now $R_\e$ and $S_\e$  can be estimated as in \cite{BNminimization} depending on whether $2^*>3$ or $2^*\leq 3$ as follows:
 \begin{equation}\label{BN02}
 R_\e, S_\e=\begin{cases}
 O(\e^{(N/2)\theta})\;\forall\; \theta<1& \text{ if }N\geq6,\\
 O(\e^2)&\text{ if } N=5,\\
 O(\e ^{2\theta})\;\forall\; \theta<1& \text{ if }N=4,\\
 O(\e)&\text{ if } N=3.
 \end{cases}
 \end{equation}
 Thus for all $N\geq 3$, we get $R_\e,S_\e=o(\e^{(N-2)/2}).$ Therefore
\begin{align}\label{2sol13}
E_\la^a(u_\la+tRU_\e)&=\frac{1}{2}\int_\Om|\nabla u_\la|^2+\frac{R^{2}t^{2}}{2}\int_\Om|\nabla U_\e|^2 +\la\int_\Om\chi_{\{u_\la<a\}}u_\la^{-\delta}tRU_\e \nonumber \\
&\quad -\la\int_\Om G(u_\la+tRU_\e)-\frac{\la}{2^*}\int_\Om u_\la^{2^*}-\frac{\la t^{2^*}R^{2^*}}{2^*}\int_\Om U_\e^{2^*}\nonumber\\
&\quad -\la R^{2^*-1}t^{2^*-1}\int_\Om  U_\e^{2^*-1}u_\la+ o_\e(\e^{(N-2)/2})\nonumber\\
&=E_\la^a(u_\la)+\frac{R^{2}t^{2}B}{2}-\frac{\la t^{2^{*}}R^{2^{*}}A}{2^{*}}
-\la R^{2^{*}-1}t^{2^{*}-1}\int_\Om U_\e^{2^*-1}u_\la\nonumber\\
&+\la\int_\Om\left(\chi_{\{u_\la<a\}}u_\la^{-\delta}tRU_\e+G(u_\la)-G(u_\la+tRU_\e)\right)+ o_\e(\e^{(N-2)/2}).
\end{align}
Now we estimate the last integral, which we denote by $T$,
on the right hand side of (\ref{2sol13}) as follows:
\begin{align}\label{singularterm}
 T&=\int_\Om\left(\chi_{\{u_\la<a\}}u_\la^{-\delta}tRU_\e+G(u_\la)-G(u_\la+tRU_\e)\right)\nonumber\\
 &\leq\int_{A_1}u_\la^{-\de}tRU_\e
 +\int_{A_2}\left(u_\la^{-\de}tRU_\e+\frac{u_\la^{1-\de}}{1-\de}-\frac{(u_\la+tRU_\e)^{1-\de}}{1-\de}\right)\nonumber
\end{align}
where $A_1=\{x\in\Om:u_\la(x)<a\leq (u_\la+tRU_\e)(x)\}$ and
$A_2=\{x\in\Om:(u_\la+tRU_\e)(x)<a\}$.
Note that as $U_\e\rar 0$ uniformly in $\{x\in\Om:|x-y|>r\}$, we get $|A_1\setminus \overline{B_r(y)}|\rar 0$
as $\e\rar 0$.  Also as $u_\la$ is continuous and $\overline{B_r(y)}\subset \Om_a :=\{x \in \Omega: u_\la<a\}$, there exists $\ga>0$ such that $u_\la<a-\ga$
in $B_r(y)$. Thus for $x \in A_1\cap B_r(y)$, $tRU_\e(x)>\ga$, i.e.,
$$\eta(x)\frac{C_N\e^{(N-2)/2}}{(\e^2+|x-y|^2)^{(N-2)/2}}\geq\frac{\ga}{tR}.$$
Therefore $|x-y|\leq \sqrt \e(\frac{tR}{\ga})^{1/(N-2)}C_N^{1/(N-2)}=r_\e.$
Thus $A_1\cap B_r(y)\subset B_{r_\e}(y)$ and
\begin{equation}\label{T1}
\int_{A_1\cap B_r(y)}U_\e\leq\int_{B_{r_\e}(y)}U_\e\leq O_\e(1)\e^{(N-2)/2}\int_0^{r_\e}rdr=O_\e(1) \e^{N/2}.
\end{equation}
If $x\in A_1\setminus B_r(y)$, then $U_\e(x)\leq O_\e(1)\frac{\e^{(N-2)/2}}{r^{N-2}}$ and hence
\begin{equation}\label{T2}
\int_{A_1\setminus B_r(y)}U_\e=O_\e(1)\frac{\e^{(N-2)/2}}{r^{N-2}}|A_1\setminus B_r(y)|.
\end{equation}
Hence from \eqref{T1} and \eqref{T2} we get
\begin{equation}\label{singterm1}
 \int_{A_1}u_\la^{-\de}tRU_\e=o(\e^{(N-2)/2}).
\end{equation}
 Also as in page 176 of \cite{GR}
\begin{equation}\label{singterm2}
 \int_{A_2}\left(u_\la^{-\de}tRU_\e+\frac{u_\la^{1-\de}}{1-\de}-\frac{(u_\la+tRU_\e)^{1-\de}}{1-\de}\right)\leq o(\e^{(N-2)/2}).
\end{equation}
Indeed, fix $0<\tau<1/4$.Then
\begin{align}\label{revised2}
 D_\e&=\int_{A_2}\left(u_\la^{-\de}tRU_\e+\frac{u_\la^{1-\de}}{1-\de}-\frac{(u_\la+tRU_\e)^{1-\de}}{1-\de}\right)\nonumber\\
 &=\int_{|x-y|\leq\e^\tau}\left(u_\la^{-\de}tRU_\e+\frac{u_\la^{1-\de}}{1-\de}-\frac{(u_\la+tRU_\e)^{1-\de}}{1-\de}\right)\nonumber\\
 &~~~+\int_{|x-y|>\e^\tau}\left(u_\la^{-\de}tRU_\e+\frac{u_\la^{1-\de}}{1-\de}-\frac{(u_\la+tRU_\e)^{1-\de}}{1-\de}\right).
\end{align}
Now for $\e$ small, we have the following estimate for the first term on the right hand side of \eqref{revised2}.
\begin{align}\label{revised3}
 &\int_{|x-y|\leq\e^\tau}\left(u_\la^{-\de}tRU_\e+\frac{u_\la^{1-\de}}{1-\de}-\frac{(u_\la+tRU_\e)^{1-\de}}{1-\de}\right)
 \nonumber\\
  &~~~\leq c_1R\int_{|x-y|\leq\e^\tau}U_\e\,dx\nonumber\\
  &~~~=c_2R\int_{|x-y|\leq\e^\tau}
  \frac{\e^{(N-2)/2}}{(\e^2+|x-y|^2)^{(N-2)/2}}dx\nonumber\\
  &~~~\leq c_2R\e^{(N-2)/2}\int_0^{\e^\tau}r\,dr\nonumber\\
  &~~~\leq c_3R\e^{(N-2)/2+2\tau}
\end{align}
Now using the fact that $u_\la$ is bounded below in the support of 
$\eta$ and the mean value theorem, we have the following estimate for the second term on the right hand side of \eqref{revised2}.
\begin{align}\label{revised4}
 &\int_{|x-y|>\e^\tau}\left(u_\la^{-\de}tRU_\e+\frac{u_\la^{1-\de}}{1-\de}-\frac{(u_\la+tRU_\e)^{1-\de}}{1-\de}\right)
 \nonumber\\
 &~~~\leq\int_{|x-y|>\e^\tau}\left(u_\la^{-\de}-(u_\la+\theta_1tRU_\e)^{-\de}\right)tRU_\e\,dx\nonumber\\
 &~~~~\leq c_4\int_{supp\,\,\eta\cap|x-y|>\e^\tau}(u_\la+\theta_2tRU_\e)^{-1-\de}(tRU_\e)^2\,dx\nonumber\\
 &~~~~\leq c_5 R^2\int_{|x-y|>\e^\tau}U_\e^2\,dx\nonumber\\
 &~~~~\leq c_6R^2\int_{|x-y|>\e^\tau}\frac{\e^{(N-2)}}{(\e^2+|x-y|^2)^{(N-2)}}dx \nonumber\\
 &~~~~\leq c_7 R^2\e^{N-2-2\tau(N-2)}
\end{align}
for some $0<\theta_1,\theta_2<1$. Thus \eqref{revised3} and \eqref{revised4} gives \eqref{singterm2}.
Thus substituting \eqref{singterm1} and \eqref{singterm2} in (\ref{2sol13}) we get
\begin{align*}
E_\la^a(u_\la+tRU_\e)&\leq E_\la^a(u_\la)+\frac{R^{2}t^{2}B}{2}-\frac{\la t^{2^{*}}R^{2^{*}}A}{2^{*}}
-\la R^{2^{*}-1}t^{2^{*}-1}\int_\Om U_\e^{2^*-1}u_\la\\
&~~~~~~~~+ o_\e(\e^{(N-2)/2}).
\end{align*}
Now the lemma follows
using the arguments of Section 3 of \cite{Tara}. \hfill{\qed}
\begin{lemma}
Let $(MP)$ hold. Then there exists a solution $v_\la\in H^+$ of $(\tilde{P}_\lambda^a)$
and hence a second solution for $(P_\la^a)$.
\end{lemma}
\noindent
\proof Here we argue as in Lemma 3.5 of \cite{BT}. Define a complete metric space $(X,d)$ as
\[X=\{\eta\in C([0,1],H^+)\, : \eta(0)=0,\|\eta(1)\|_{\HH}>\rho_1,\, I_\lambda(\eta(1))<0\},\]
\begin{equation*}
d(\eta,\chi)=\displaystyle\max_{t\in[0,1]}\;\|\eta(t)-\chi(t)\|_{\HH}.
 \end{equation*}
From $(i)$ of Lemma (\ref{L2}), if $R$ is chosen large, it is clear that $X$ is non-empty.
Let
$\gamma_0=\displaystyle\inf_{\eta\in X}\;\displaystyle\max_{s\in[0,1]}\;I_\lambda(\eta(s)).$
Then $(ii)$ of Lemma (\ref{L2}) and $(MP)$ implies that
\begin{equation}
\label{energy level}0<\ga_0<\displaystyle\frac{S^{\frac{N}{2}}}{N\la^{(N-2)/2}}.
\end{equation}
 Define
\begin{equation*}\label{eq3.20}
\Psi(\eta)=\displaystyle\max_{t\in[0,1]}\;I_\lambda(\eta(t)),\eta\in X.
 \end{equation*}
Thus applying Ekeland's variational principle to the above functional
we get a sequence $\{\eta_n\}\subseteq X$ such that
\begin{itemize}
 \item [(i)]$\displaystyle\max_{t\in[0,1]}I_\la(\eta_n(t))<\ga_0+\frac{1}{n}$.
\item[(ii)]$\displaystyle\max_{t\in[0,1]}I_\la(\eta_n(t))\leq\displaystyle\max_{t\in[0,1]}I_\la(\eta(t))
+\frac{1}{n}\displaystyle\max_{t\in[0,1]}\|\eta(t)-\eta_n(t)\|_{\HH}$ for all $\eta\in X$.
\end{itemize}
Set $\Lambda_n=\{t\in[0,1]|I_\la(\eta_n(t))=\displaystyle\max_{s\in[0,1]}I_\la(\eta_n(s))\}$.
Then as in the Claim on page 659 of \cite{BT}
we get $t_n\in\La_n$ such that for $v_n=\eta_n(t_n)$ and $\xi\in H^+$ we have
\begin{equation}\label{eqn 2.26}
I_\la^0\left(v_n, \frac{\xi-v_n}{max\{1,\|\xi-v_n\|_{\HH}\}}\right)\geq-\frac{1}{n}
\end{equation}
and
\begin{equation}\label{eqn 2.27}
I_\la(v_n)\rightarrow \ga_0 ~~\text{ as } n\rightarrow\infty.
\end{equation}
From (\ref{eqn 2.27}) we have
$$\displaystyle\frac{1}{2}\|v_n\|_{\HH}^2-\la\int_\Om \tilde F(v_n)-\la\int_\Om \tilde G(v_n)\leq\ga_0+o_n(1).$$
As $\tilde G(v_n)\leq0$ and $\tilde F(v_n)=\displaystyle\frac{(v_n+u_\la)^{2^*}}{2^*}-u_\la^{2^*-1}v_n-\frac{u_\la^{2^*}}{2^*}$, this implies
\begin{equation}\label{eqn 2.28}
\frac{1}{2}\|v_n\|_{\HH}^2-\la\int_\Om \frac{(v_n+u_\la)^{2^*}}{2^*}+\la\int_\Om u_\la^{2^*-1}v_n\leq\ga_0+o_n(1).
\end{equation}
Also substituting  $\xi=2v_n+u_\la$ in (\ref{eqn 2.26}), by Remark \ref{w-existence} we obtain  (by abuse of notation) 
$$w^{v_n}_n(x)\in[\chi_{\{v_n(x)+u_\la(x)<a\}},\chi_{\{v_n(x)+u_\la(x)\leq a\}}]$$ such that
\begin{equation}\label{eqn 2.29}
\|v_n+u_\la\|_{\HH}^2-\lambda\int_{\Omega}((v_n+u_\la)^{2^*}-
 w^{v_n}_n( v_n+u_\la)^{1-\de})\geq-\frac{1}{n}\max\;\{1,\|v_n+u_\la\|_{\HH}\}.
 \end{equation}
From (\ref{eqn 2.28}) and (\ref{eqn 2.29}) we derive
\begin{equation*}\label{eqn 2.30}
 \frac{1}{2}\|v_n\|_{\HH}^2-\frac{1}{2^*}\|v_n\|_{\HH}^2\leq c_1+c_2\|v_n\|_{\HH}
\end{equation*}
where $c_1,c_2>0$. Thus $\|v_n\|_{\HH}\leq C$ for all $n\in\mathbb N$.
Hence $v_n\rightharpoonup v_\la$ weakly in $H^1_0(\Om)$ and as in case of zero altitude $v_\la$ solves $(\tilde P_\la^a)$.
Now we claim that $v_n\rightarrow v_\la$ in $H^1_0(\Om)$ and thus $v_\la\not\equiv0$.
{Without loss of generality we assume $I_\la(v_\la)=0,$ otherwise it would
imply that $v_\la\not \equiv 0$ and we are done.}
As $\|v_n\|_{\HH}\leq C$, from (\ref{eqn 2.26}), for $\xi\in H^+$
we have $I_\la^0(v_n, \xi-v_k)\geq -\displaystyle\frac{C_1}{n}(1+\|\xi\|_{\HH})=o_n(1)$.
Then as in zero altitude case we get
\begin{equation}\label{eqn 2.31}
 \|v_n-v_\la\|_{\HH}^2-\la\|v_n-v_\la\|^{2^*}_{L^{2^*}(\Om)}=o_n(1).
\end{equation}
Also by Brezis-Lieb lemma,
\begin{align}\label{eqn 2.32}
 I_\la(v_n)&=\frac{1}{2}\|v_n\|_{\HH}^2-\la\int_\Om\tilde F(v_n)-\la\int_\Om\tilde G(v_n)\nonumber\\
&=\frac{1}{2}\|v_n-v_\la\|_{\HH}^2+\frac{1}{2}\|v_\la\|_{\HH}^2+\int_\Om\nabla(v_n-v_\la)\cdot\nabla v_\la\nonumber\\
&~~~~-\la\left(\frac{1}{2^*}\int_\Om(v_n+u_\la)^{2^*}-\frac{1}{2^*}\int_\Om u_\la^{2^*}-\int_\Om u_\la^{2^*-1}v_n\right)
-\la\int_\Om\tilde G(v_n)\nonumber\\
&=\frac{1}{2}\|v_n-v_\la\|_{\HH}^2-\frac{\la}{2^*}\|v_n-v_\la\|_{L^{2^*}(\Om)}^{2^*}
+\frac{1}{2}\|v_\la\|_{\HH}^2-\la\int_\Om\tilde G(v_n)\nonumber\\
 &~~~~-\la\left(\frac{1}{2^*}\int_\Om(v_\la+u_\la)^{2^*}-\frac{1}{2^*}\int_\Om u_\la^{2^*}-\int_\Om u_\la^{2^*-1}v_n\right)\nonumber\\
&~~~~+\int_\Om\nabla(v_n-v_\la)\cdot\nabla v_\la+o_n(1)\nonumber\\
&=\frac{1}{2}\|v_n-v_\la\|_{\HH}^2-\frac{\la}{2^*}\|v_n-v_\la\|_{L^{2^*}(\Om)}^{2^*}
+I_\la(v_\la)\nonumber\\
&~~~~+\la\int_\Om\left(\tilde G(v_\la)-\tilde G(v_n)\right)+o_n(1)\nonumber\\
&=\frac{1}{2}\|v_n-v_\la\|_{\HH}^2-\frac{\la}{2^*}\|v_n-v_\la\|_{L^{2^*}(\Om)}^{2^*}
+I_\la(v_\la)+o_n(1).
\end{align}
Now as $I_\la(v_\la)= 0$, using (\ref{energy level}),
(\ref{eqn 2.27}), (\ref{eqn 2.31}) and (\ref{eqn 2.32}), we get
\begin{align}\label{eqn 2.33}
 \|v_n-v_\la\|^2_{\HH}=N\gamma_0+o_n(1)<\frac{S^{\frac{N}{2}}}{\la^{(N-2)/2}}+o_n(1).
\end{align}
Also by the Sobolev embedding we have
\begin{align}\label{eqn 2.34}
 \|v_n-v_\la\|^2_{\HH}\left(1-\la S^{\frac{-2^*}{2}}\|v_n-v_\la\|_{\HH}^{2^*-2}\right)
&\leq \|v_n-v_\la\|_{\HH}^2-\la\|v_n-v_\la\|_{L^{2^*}(\Om)}^{2^*}\nonumber\\
&=o_n(1).
\end{align}
Thus combining (\ref{eqn 2.33}) and (\ref{eqn 2.34}) we obtain
 $v_n\rar v_\la$ in $H^1_0(\Om)$.
This completes the proof of the lemma.\hfill{\qed}\\
We are now ready to give the

{\bf Proof of Theorem 1.1:} The existence of the first solution $u_\la$ for all $\la \in (0,\La^a)$ 
follows from lemma 2.1 and theorem 2.2 . The existence of the second solution $v_\la$ for 
the same range of $\la$ follows from lemma 3.1 and lemma 3.3 keeping in view the remark 3.3. \hfill \qed

{\bf Acknowledgement:} We would like to thank the anonymous referee for his meticulous review which greatly improved the presentation of the paper.

\end{document}